\documentclass[11pt]{amsart}
\usepackage{amsmath,amssymb}

\begin{document}

\newtheorem{thm}{Theorem}[section]
\newtheorem{lem}[thm]{Lemma}
\newtheorem{prop}[thm]{Proposition}
\newtheorem{cor}[thm]{Corollary}
\newtheorem{defn}[thm]{Definition}
\newtheorem*{remark}{Remark}

\numberwithin{equation}{section}

\newcommand{\Z}{{\mathbb Z}} 
\newcommand{\Q}{{\mathbb Q}}
\newcommand{\R}{{\mathbb R}}
\newcommand{\C}{{\mathbb C}}
\newcommand{\N}{{\mathbb N}}
\newcommand{\FF}{{\mathbb F}}
\newcommand{\fq}{\mathbb{F}_q}

\newcommand{\rmk}[1]{\footnote{{\bf Comment:} #1}}

\renewcommand{\mod}{\;\operatorname{mod}}
\newcommand{\ord}{\operatorname{ord}}
\newcommand{\TT}{\mathbb{T}}
\renewcommand{\i}{{\mathrm{i}}}
\renewcommand{\d}{{\mathrm{d}}}
\renewcommand{\^}{\widehat}
\newcommand{\HH}{\mathbb H}
\newcommand{\Vol}{\operatorname{vol}}
\newcommand{\area}{\operatorname{area}}
\newcommand{\tr}{\operatorname{tr}}
\newcommand{\norm}{\mathcal N} 
\newcommand{\intinf}{\int_{-\infty}^\infty}
\newcommand{\ave}[1]{\left\langle#1\right\rangle} 
\newcommand{\Var}{\operatorname{Var}}
\newcommand{\Prob}{\operatorname{Prob}}
\newcommand{\sym}{\operatorname{Sym}}
\newcommand{\disc}{\operatorname{disc}}
\newcommand{\CA}{{\mathcal C}_A}
\newcommand{\cond}{\operatorname{cond}} 
\newcommand{\lcm}{\operatorname{lcm}}
\newcommand{\Kl}{\operatorname{Kl}} 
\newcommand{\leg}[2]{\left( \frac{#1}{#2} \right)}  

\newcommand{\sumstar}{\sideset \and^{*} \to \sum}

\newcommand{\LL}{\mathcal L} 
\newcommand{\sumf}{\sum^\flat}
\newcommand{\Hgev}{\mathcal H_{2g+2,q}}
\newcommand{\USp}{\operatorname{USp}}

\title[Statistics for zeros of hyperelliptic zeta functions]
{Statistics of the zeros of zeta functions in families of
hyperelliptic curves over a finite field}
 \author{Dmitry Faifman}
 \email{faifmand@post.tau.ac.il}

\author{Ze\'ev Rudnick}
\email{rudnick@post.tau.ac.il}

\address{Raymond and Beverly Sackler School of Mathematical Sciences,
Tel Aviv University, Tel Aviv 69978, Israel}

\date{May 4, 2008}

 \thanks{Supported by the Israel Science Foundation (grant No. 925/06).}

\begin{abstract}
We study the fluctuations in the distribution of zeros of zeta
functions of a family of 
hyperelliptic curves defined over a fixed finite field, in the
limit of large genus.  
According to the Riemann Hypothesis for curves, the zeros all lie on a
circle. Their angles are uniformly distributed, so for a curve of
genus $g$ a fixed interval $\mathcal I$ will contain asymptotically 
$2g|\mathcal I|$
angles as the genus grows. We show that for the variance of number of 
angles in $\mathcal I$ is asymptotically $\frac
2{\pi^2}\log(2g|\mathcal I|)$ and prove a central limit theorem: The
normalized fluctuations are Gaussian. These results continue to hold
for shrinking intervals as long as the expected number of angles
$2g|\mathcal I|$ tends to infinity. 

\end{abstract}

\maketitle

\section{Introduction}

 Let $C$ be a smooth, projective, geometrically connected curve of genus
$g\geq 1$  defined over a  finite field $\FF_q$ of cardinality  $q$.
The zeta function of the curve is defined as
\begin{equation}\label{zeta function}
Z_C(u):=\exp \sum_{n=1}^\infty N_n\frac{u^n}{n}, \quad |u|<1/q
\end{equation}
where $N_n$ is the number of points on $C$ with coefficients in an
extension $\FF_{q^n}$ of $\FF_q$ of degree $n$.
The zeta function is a rational function of the form
$$Z_C(u) = \frac{P_C(u)}{(1-u)(1-qu)}$$
where $P_C(u)\in \Z[u]$ is a polynomial of degree $2g$, with $P(0)=1$,
satisfies the functional equation
$$ P_C(u) = (qu^2 )^{g} P_C(\frac 1{qu})$$
and has all its zeros on the circle $|u|=1/\sqrt{q}$ (this is the
Riemann Hypothesis for curves \cite{Weil}).
Moreover, there is a unitary  symplectic
matrix $\Theta_C\in \USp(2g)$, defined up to conjugacy, so that
$$P_C(u) = \det (I-u\sqrt{q} \Theta_C)$$
The eigenvalues of $\Theta_C$ are of the form
$e^{2\pi i \theta_{C,j}}$, $j=1,\dots, 2g$.

Our goal is to study the statistics of the set of angles
$\{\theta_{j,C}\}$ as we draw $C$ at random from a family
of hyperelliptic curves of genus $g$ defined over $\FF_q$ where
$q$ is assumed to be odd.
The family, denoted by
$\Hgev$,  is that of curves having an affine equation of
the form $y^2=Q(x)$, with $Q\in \FF_q[x]$ a monic,  square-free
polynomial of degree $2g+2$. The corresponding function field is
called a real quadratic  function field.
The measure on $\Hgev$ is simply the uniform probability
measure on the set of such polynomials $Q$.

A fundamental statistic is the counting function of the angles. Thus
for an interval\footnote{Due to the functional equation, it
suffices to restrict the discussion to symmetric intervals.}
 $\mathcal I = [-\frac \beta 2,\frac \beta 2]$ (which may vary
with the genus $g$ or with $q$), let
$$N_{\mathcal I}(C) = \#\{j: \theta_{j,C} \in \mathcal I \}$$
The angles are uniformly distributed as $g\to \infty$
(see Proposition~\ref{prop:unif dist}): For fixed $\mathcal I$,
$$
N_{\mathcal I}(C)\sim 2g|\mathcal I| \;.
$$
We wish to study the fluctuations of $N_{\mathcal I}$ as we vary $C$ in
$\Hgev$.
This is in analogy to the work of Selberg \cite{Se1, Se2, Se3}, who
studied the  fluctuations
in the number $N(t)$ of zeros of the Riemann  zeta function $\zeta(s)$
up to height $t$. By the Riemann-von Mangoldt formula,
$$ N(t) =  \frac{t}{2\pi}\log \frac{t}{2\pi e} + \frac 78
+S(t)+O(\frac 1t)
$$
with $S(t) = \frac 1\pi \arg \zeta(\frac 12+it)$. Selberg showed that
the variance of $S(t)$, for $t$ picked uniformly in $[0,T]$, is
$\frac 1{2\pi^2} \log\log T$, and that the moments of
$S(t)/\sqrt{\frac 1{2\pi^2}\log\log t}$ are those of a standard Gaussian.

Katz and Sarnak \cite{KS} showed that for fixed genus,
the conjugacy classes $\{\Theta_C: C\in \Hgev \}$ become
uniformly distributed in $\USp(2g)$ in the limit  $q\to\infty$
of large constant field size.
In particular the statistics of $N_{\mathcal I}$ are the same
as those of the corresponding quantity for a random matrix in
$\USp(2g)$. That is, if $U\in \USp(2g)$ is a unitary symplectic
matrix, with eigenvalues $e^{2\pi i \theta_j(U)}$, $j=1,\dots ,2g$,
set
$$ \^N_{\mathcal I}(U) = \#\{j: \theta_j(U) \in \mathcal I \}$$
Then the work of Katz and Sarnak \cite{KS} gives
\begin{equation}\label{eq:KS}
\lim_{q\to \infty}
\Prob_{\Hgev} \left(N_{\mathcal I}(C)=k \right)
=
\Prob_{\USp(2g)} \left( \^N_{\mathcal I}(U) =k \right)
\end{equation}

In the limit of large matrix size, the statistics of 
$\^N_{\mathcal I}(U)$ and related quantities, such 
as the logarithm of the characteristic polynomial of $U$, 
have been found to have Gaussian fluctuations  in various ensembles of
random matrices \cite{Politzer, 
Costin-Lebowitz, Baker-Forrester, Johansson, Keating-Snaith, Soshnikov, DE,  
HKO, Wieand}. In particular, when averaged over
$\USp(2g)$, the expected value of $\^N_I$ is 
$2g|\mathcal I|$,  
the variance is $\frac 2{\pi^2}\log(2g|\mathcal I|)$ and the
normalized random variable $(\^N_{\mathcal I} - 2g|\mathcal I|) 
/\sqrt{\frac 2{\pi^2}\log(2g|\mathcal I|)}$ has a normal distribution
as $g\to \infty$.   
Moreover this holds for shrinking intervals, that is
if we take the length of the interval $|\mathcal I|\to 0$ as $g\to \infty$ as
long as the expected number of angles tends to
infinity\footnote{This is sometime called the ``mesoscopic'' regime},
that is as long as $2g|\mathcal I|\to \infty$. 
Thus \eqref{eq:KS} implies that for the iterated limit
$\lim_{g\to\infty}(\lim_{q\to\infty})$ we get a Gaussian distribution:
\begin{equation*}
\lim_{g\to\infty} \left( \lim_{q\to\infty}
\Prob_{\Hgev} \left(a< \frac{N_{\mathcal I}(C)-2g|\mathcal I|}
{\sqrt{\frac 2{\pi^2} \log(2g|\mathcal I|)}}<b \right) \right)
 = \frac 1{\sqrt{2\pi}} \int_a^b e^{-x^2/2}dx
\end{equation*}

In this paper we will study these problems for a {\em fixed} constant
field $\FF_q$ in the limit of large genus
$g\to \infty$, that is without first taking $q\to\infty$,
which was crucial to the approach of Katz and Sarnak.
We will show that as $g\to \infty$, for both the global regime
($|\mathcal I|$ fixed) and the mesoscopic regime ($|\mathcal I|\to 0$
while $2g|\mathcal I|\to \infty$), the expected value of
$N_{\mathcal I}$ is $2g|\mathcal I|$,
the variance  is asymptotically
$\frac 2{\pi^2}\log(2g|\mathcal I| ) $
and that the fluctuations are Gaussian, that is for fixed $a<b$,
\begin{equation}\label{CLT for S}
\lim_{g\to\infty}
\Prob_{\Hgev}\left(a< \frac{N_{\mathcal I}-2g|\mathcal I|}{\sqrt{\frac 2{\pi^2} \log(
    2g|\mathcal I|)}}<b \right)
= \frac 1{\sqrt{2\pi}} \int_a^b e^{-x^2/2}dx
\end{equation}

Our argument hinges upon the fact that $P_C(u)$ is
the L-function attached to a quadratic character of
$\FF_q[x]$.  Thus for $Q$ monic, square free, of degree $2g+2$ the
quadratic character $\chi_Q$ is defined in terms of the quadratic
residue symbol as $\chi_Q(f) = \leg{Q}{f}$
(see \S~\ref{subsec:quad}). The associated L-function is
$$ \LL(u,\chi_Q) = \prod_P (1-\chi_Q(u)u^{\deg P})^{-1}$$
the product taken over all monic irreducible polynomials $P\in
\FF_q[x]$.
Then
$$P_C(u)=(1-u)^{-1} \LL(u,\chi_Q)$$
as was found in E. Artin's thesis \cite{Artin}.
Thus one may tackle the problem using Selberg's
original arguments \cite{Se1}\footnote{The paper \cite{Se1} is under the
Riemann hypothesis;  \cite{Se2, Se3} are unconditional.}
adapted to the function field setting; this was carried out in the
M.Sc. thesis of the first-named author \cite{Faifman}.
Instead we follow a quicker route, via the explicit formula, used
recently by Hughes, Ng and Soundararajan \cite{HNS}.

An important challenge is to investigate the local regime, when the
length of the interval is of order $1/2g$ as $g\to \infty$. 
Due to the Central Limit Theorem for random matrices, 
we may rewrite
\eqref{CLT for S} as
\begin{multline}\label{equivalent form of CLT}
\lim_{g\to\infty}
\Prob_{\Hgev}\left(a< \frac{N_{\mathcal I}-2g|\mathcal I|}
{\sqrt{\frac 2{\pi^2} \log( 2g|\mathcal I|)}}<b \right)
\\ = \lim_{g\to\infty} \Prob_{\USp(2g)}
\left(a< \frac{\^N_{\mathcal I}-2g|\mathcal I|}
{\sqrt{\frac 2{\pi^2} \log( 2g|\mathcal I|)}}<b \right)
\end{multline}
and ask if \eqref{equivalent form of CLT} remains valid also for
shrinking intervals of the form $\mathcal I=\frac 1{2g}\mathcal J$ where
$\mathcal J$ is fixed, when the result is no longer a Gaussian.
An equivalent form of \eqref{equivalent form of CLT} was conjectured in
\cite{KS-BAMS}.



\noindent{\bf Acknowledgement:} We thank Chris Hughes, Jon Keating,  
Emmanuel Kowalski and Igor Shparlinski for discussions and comments on
earlier versions of the paper. 

\section{Background on Dirichlet characters and L-functions}

\subsection{}
We review some generalities about Dirichlet L-functions for the
rational function field; see \cite{Rosen} for details.

The norm of a nonzero polynomial $f\in \FF_q[x]$ is defined as
$||f|| =q^{\deg f}$.
The zeta function of the rational function field
is
$$ \zeta_q(s):= \prod_P (1-||P||^{-s})^{-1},\quad \Re(s)>1 $$
the product over all irreducible monic polynomials (``primes'') in
$\FF_q[x]$. In terms of the more convenient variable
$$u=q^{-s}$$
the zeta function becomes
$$
Z(u) = \prod_P (1-u^{\deg P})^{-1}, \quad |u|<1/q \;.
$$
By the fundamental theorem of arithmetic in $\FF_q[x]$, $Z(u)$
can be expressed as a sum over all monic polynomials:
$$ Z(u) = \sum_{f \mbox{ monic}} u^{\deg f } $$
and hence
$$ Z(u) = \frac 1{1-qu} \;.$$

Given a monic polynomial $Q\in \FF_q[x]$, a Dirichlet character
modulo $Q$ is a homomorphism
$$\chi: (\FF_q[x]/ Q\FF_q[x])^\times \to \C^\times$$
A character modulo $Q$ is {\em primitive} if there is no proper
divisor $\tilde Q$ of $Q$ and some character $\tilde \chi$ mod $\tilde
Q$ so that $\chi(n)=\tilde\chi(n)$ whenever $\gcd(n,Q)=1$.

For a Dirichlet character $\chi$ modulo $Q$ of $\FF_q[x]$, we form the
L-function
\begin{equation}\label{euler prod}
\LL(u,\chi) = \prod_P (1-\chi(P) u^{\deg P})^{-1}
\end{equation}
(convergent for $|u|<1/q$), where $P$ runs over all monic
irreducible polynomials.
It can be expressed as a series
\begin{equation} \label{series for L}
 \LL(u,\chi)  = \sum_f \chi(f) u^{\deg f}
\end{equation}
where the sum is over all monic polynomials. If $\chi$ is nontrivial,
then it is easy to show that
$$\sum_{\deg f=n} \chi(f) = 0,\quad n\geq \deg Q$$
and hence the L-function is in fact a polynomial of degree at most $\deg Q-1$.

One needs to
distinguish ``even'' characters from the rest, where ``even''
means $\chi(cH) = \chi(H)$, $\forall c\in \FF_q^\times$. The
analogue for ordinary Dirichlet characters is $\chi(-1)=1$. For
even characters, the L-function has a trivial zero at $u=1$.

We assume from now on that $\deg Q>0$ and that $\chi$ is primitive.
One then defines a ``completed'' L-function
$$\LL^*(u,\chi) = (1-\lambda_\infty(\chi) u)^{-1} \LL(u,\chi)$$
where $\lambda_\infty(\chi)=1$ if $\chi$ is ``even'', and is zero
otherwise.
The completed L-function $\LL^*(u,\chi)$ is then a polynomial of
degree
$$D=\deg Q-1-\lambda_\infty(\chi)$$
and satisfies the functional equation
$$\LL^*(u,\chi) = \epsilon(\chi) (q^{1/2}u)^D
\LL^*(\frac 1{qu},\chi^{-1})
$$
with $|\epsilon(\chi)|=1$.
We express  $\LL^*(u,\chi)$ in term of its inverse zeros
as
\begin{equation}\label{prod zeros}
\LL^*(u,\chi) = \prod_{j=1}^D (1-\alpha_{j,\chi} u) \;.
\end{equation}
The Riemann Hypothesis in this setting, proved by Weil \cite{Weil},
is that all $|\alpha_{j,\chi}|=\sqrt{q}$.
We may thus write
\begin{equation}\label{angles of L}
\alpha_{j,\chi} = \sqrt{q}e^{2\pi i\theta_{j,\chi}} 
\end{equation}
for suitable phases $\theta_{j,\chi}\in \R/\Z$.
As a consequence, for any nontrivial character, not necessarily
primitive,  the inverse zeros of
the L-function all have absolute value $\sqrt{q}$ or $1$.

\begin{lem}\label{lem:PV}
  Let $\chi$ be a non-trivial Dirichlet character modulo $f$. Then for
  $n<\deg f$,
$$ \left| \sum_{\deg B=n} \chi(B) \right| \leq \binom{\deg f-1}{n} q^{n/2}$$
(the sum over all monic polynomials of degree $n$).
\end{lem}
\begin{proof}
Indeed, all we need to do is compare the series expansion
 \eqref{series for L} of
 $\LL(u,\chi)$, which is a polynomial of degree at most $\deg f-1$,
 with the expression in terms of the inverse zeros:
$$
\sum_{0\leq n<\deg f} (\sum_{\deg B=n} \chi(B))u^n  =
\prod_{j=1}^{\deg f-1} (1-\alpha_j u)
$$
to get
$$ \sum_{\deg B=n} \chi(B) =(-1)^n \sum_{\substack{S\subset \{1,\dots, \deg
f-1\} \\ \#S=n }} \prod_{j\in S} \alpha_j
$$
and then use $|\alpha_j|\leq \sqrt{q}$.
\end{proof}
Note that for $n\geq \deg f$ the character sum vanishes.


\subsection{Quadratic characters} \label{subsec:quad}
We assume from now on  that $q$ is odd.
Let $P(x)\in \FF_q[x]$ be  monic and irreducible. The quadratic residue symbol
$\leg{f}{P} \in \{\pm 1\}$  is
defined for $f$ coprime to $P$ by
$$\leg{f}{P} \equiv f^{\frac{||P||-1}2} \mod P \;.
$$
For arbitrary monic $Q$, the Jacobi symbol $\leg{f}{Q}$ is defined for
$f$ coprime to $Q$ by writing $Q=\prod P_j$ as a product of monic
irreducibles and setting
$$\leg{f}{Q} = \prod \leg{f}{P_j} \;.
$$
If $f,Q$ are not coprime we set $\leg{f}{Q}=0$.
If $c\in \FF_q^\ast$ is a scalar then
\begin{equation}\label{scalars}
\leg{c}{Q} = c^{\frac{q-1}2 \deg Q} \;.
\end{equation}
The law of quadratic reciprocity asserts that if $A, B\in \FF_q[x]$ are monic
and coprime then
\begin{equation}\label{quad reciprocity}
\leg{A}{B} =\leg{B}{A}(-1)^{\frac{q-1}2 \deg A\deg B}
= \leg{B}{A}(-1)^{\frac{||A||-1}2\cdot \frac{||B||-1}2}
 \;.
\end{equation}
This relation continues to hold if $A$ and $B$ are not coprime as both
sides vanish.

Given a square-free $Q\in \FF_q[x]$, we define the quadratic character
$\chi_Q$ by
$$ \chi_Q(f) = \leg{Q}{f}$$
If $\deg Q$ is even, this is a primitive Dirichlet character modulo $Q$.
Note that by virtue of \eqref{scalars}, $\chi_Q$ is an even character
(that is trivial on scalars) if and only if $\deg Q$ is even.

It is important for us that the numerator $P_C(u)$ of the zeta
function \eqref{zeta function} of the hyperelliptic curve $y^2=Q(x)$
coincides with the completed Dirichlet L-function $\LL^*(u,\chi_Q)$
associated with the quadratic character $\chi_Q$.

\subsection{The Explicit Formula}
\begin{lem}
Let $h(\theta)=\sum_{|k|\leq K} \^h(k)e(k\theta)$ be a
trigonometric polynomial, which we assume is real valued and even:
$h(-\theta)=h(\theta) = \overline{h(\theta)}$. Then for a primitive
character $\chi$ we have
\begin{multline}\label{explicit formula}
\sum_{j=1}^D h(\theta_{j,\chi})   =D\int_0^1 h(\theta)d\theta  +
\lambda_\infty(\chi) \frac 1{\pi i} \int_0^1 h(\theta)
\frac{d}{d\theta} \log (1-\frac{e^{2\pi i\theta}}{\sqrt{q}}) d\theta
\\-\sum_f \^h(\deg f) \frac{\Lambda(f)}{||f||^{1/2}}
\left(\chi(f) +\overline{\chi(f)} \right)
\end{multline}
\end{lem}
\begin{proof}
By computing the logarithmic derivative $u\frac{\LL'}{\LL}$ in two
different ways, either using the Euler product \eqref{euler prod} or
the zeros \eqref{prod zeros} we get an identity, for $n>0$,
\begin{equation*}
-\sum_{j=1}^D \alpha_{j,\chi}^n = \sum_{\deg f=n}\Lambda(f)\chi(f)
+\lambda_\infty(\chi)
\end{equation*}
where $\Lambda(f) = \deg P$ if $f=P^k$ is a prime power, and
$\Lambda(f)=0$ otherwise.
Therefore we get an explicit formula  in terms of the phases $\theta_{j,\chi}$
\begin{equation*}
-\sum_{j=1}^D e^{2\pi i n\theta_{j,\chi}} =
\frac{\lambda_\infty(\chi)}{q^{|n|/2}} +
\sum_{\deg f=|n|}\frac{\Lambda(f)}{||f||^{1/2}}
\begin{cases}   \overline{\chi(f)}&n<0\\ \chi(f)& n>0    \end{cases}
\end{equation*}
which is valid  for $n$ both positive and negative.

Now let $h(\theta)=\sum_{|k|\leq K} \^h(k)e(k\theta)$ be a
trigonometric polynomial, which we assume is real valued and even:
$h(-\theta)=h(\theta) = \overline{h(\theta)}$. Then the Fourier
coefficients are also real and even:
$\^h(-k)=\^h(k)=\overline{\^h(k)}$.
Using the Fourier expansion of $h$ we get
\begin{multline*}
\sum_{j=1}^D h(\theta_j)  = D\^h(0) +\sum_j \sum_{k=1}^K \^h(k)
(e(k\theta_j) +e(-k\theta_j)) \\
=  D\int_0^1 h(\theta)d\theta  - \sum_{k=1}^K \^h(k) \left(
2\frac{\lambda_\infty(\chi)}{q^{k/2}} +
\sum_{\deg f=k} \frac{\Lambda(f)}{||f||^{1/2}}
\left(\chi(f) +\overline{\chi(f)} \right) \right) \\
=D\int_0^1 h(\theta)d\theta  -
2\lambda_\infty(\chi) \sum_{k=1}^K \frac{\^h(k)}{q^{k/2}}
-\sum_f \^h(\deg f) \frac{\Lambda(f)}{||f||^{1/2}}
\left(\chi(f) +\overline{\chi(f)} \right)
\end{multline*}
Note that since $h$ is real valued,
$$\sum_{k=1}^K \frac{\^h(k)}{q^{k/2}}  = \int_0^1 h(\theta)
\frac{q^{-1/2}e^{2\pi i\theta}}{1 - q^{-1/2}e^{2\pi i\theta}} =  \frac
1{2\pi i} \int_0^1 h(\theta)
\frac{d}{d\theta} \log \frac 1{1-\frac{e^{2\pi i\theta}}{\sqrt{q}}}
d\theta
$$
which gives the claim.
\end{proof}
\subsection{}
For the quadratic character $\chi_Q$, with $Q$ square-free of
degree $2g+2$, we get $\lambda_\infty=1$, $D=2g$,  and
the explicit formula reads
\begin{multline}\label{explicit formula quadratic}
\sum_{j=1}^{2g} h(\theta_{j,Q})   =2g \int_0^1 h(\theta)d\theta
+\frac 1{\pi i} \int_0^1 h(\theta)
\frac{d}{d\theta} \log (1-\frac{e^{2\pi i\theta}}{\sqrt{q}}) d\theta
\\
-2\sum_f \^h(\deg f) \frac{\Lambda(f)}{||f||^{1/2}} \chi_Q(f)
\end{multline}

\section{Averaging over $\Hgev$}
Let $\mathcal H_{d,q}\subset \FF_q[x]$ be the set of all square-free
monic polynomials of degree $d$. The cardinality of $\mathcal H_{d,q}$ is
\begin{equation*}
\# \mathcal H_{d,q} = \begin{cases}
(1-\frac 1q)q^{d},& d\geq 2\\q,&d=1
\end{cases}
\end{equation*}
as may be seen by expressing  the generating function
$\sum_{d=0}^\infty \mathcal H_{d,q} u^d$ in terms of the zeta function
$Z(u)$ of the rational function field:
$$ Z(u)=Z(u^2) \sum_{d=0}^\infty \mathcal H_{d,q} u^d
$$
In particular we have
\begin{equation}\label{card Hgev}
\#\Hgev = (1-\frac 1q) q^{2g+2}
\end{equation}

We denote by $\ave{\bullet}$ the mean value of any quantity defined on
$\Hgev$, that is
$$\ave{ F}:= \frac 1{\#\Hgev}
\sum_{Q\in \Hgev} F(Q)
$$
\begin{lem}\label{lem:average chi}
If $f\in \FF_q[x]$ is not a square then
$$ \ave{\chi_Q(f)} \leq \frac{2^{\deg f-1}}{(1-\frac 1q) q^{g+1}}$$
\end{lem}
\begin{proof}
We use the Mobius function to pick out the square free monic
polynomials via the formula
$$ \sum_{A^2\mid Q}\mu(A) =
\begin{cases}  1,& Q
  \mbox{ square-free}\\0,&\mbox{otherwise} \end{cases}
$$
where we sum over all monic polynomials whose square divides $Q$.
Thus the sum over all square-free polynomials is given by
\begin{equation*}  \begin{split}
\sum_{Q\in \Hgev} \chi_Q(f) &= \sum_{\deg Q=2g+2} \sum_{A^2\mid
  Q} \mu(A)\leg{Q}{f} \\
&= \sum_{\deg A\leq g+1} \mu(A) \leg{A}{f}^2 \sum_{\deg B=2g+2-2\deg A}
\leg{B}{f}
  \end{split}
\end{equation*}
To deal with the inner sum,
note that $\leg{\bullet}{f}$ is a non-trivial character since $f$ is not a
square,  so we can use
Lemma~\ref{lem:PV} to get
\begin{equation}
  \left|\sum_{\deg B=2g+2-2\deg A}\leg{B}{f} \right|
\leq \binom{\deg f-1}{2g+2-2\deg A} q^{g+1 -\deg A}
\end{equation}
if $2g+2-2\deg A <\deg f$, and the sum is zero otherwise.
Hence we have
\begin{equation*}
  \begin{split}
\left| \sum_{Q\in \Hgev} \chi_Q(f) \right| &\leq
\sum_{\deg A\leq g+1} \left|\sum_{\deg B=2g+2-2\deg A}\leg{B}{f} \right| \\
& \leq \sum_{g+1-\frac{\deg f}2< \deg A\leq g+1}
\binom{\deg f-1}{2g+2-2\deg A} q^{g+1 -\deg A}\\
& =q^{g+1 } \sum_{ g+1-\frac{\deg f}2< j \leq g+1}\binom{\deg f-1}{2g+2-2j}
\leq 2^{\deg f-1}q^{g+1 }
      \end{split}
\end{equation*}
Dividing by $\#\Hgev = q^{2g+2}(1-\frac 1q)$ proves the lemma.
\end{proof}

\begin{lem}\label{lem:ave chi^2}
Let $P_1,...,P_k$ be prime polynomials. Then
$$
\ave{\chi_Q(\prod_{j=1}^k P_j^2)}= 
1 +O\left(\sum_{j=1}^k\frac 1{||P_j||}\right) 
$$
\end{lem}
\begin{proof}
We have $\chi_Q(\prod_{j=1}^k P_j^2)=1$ if
$\gcd(\prod_{j=1}^k P_j,Q)=1$, and $\chi_Q(\prod_{j=1}^k P_j^2)=0$
otherwise. Since for primes $P_1,...,P_k$ the condition
$\gcd(\prod_{j=1}^k P_j,Q)\neq 1$ is equivalent to $P_j$ dividing
$Q$ for some $j$, we may write
$$\chi_Q(\prod_{j=1}^k P_j^2) =1-  \begin{cases}
1,&\exists P_j\mid Q\\0,& \mbox{otherwise}
\end{cases}
$$
and hence
$$\ave{\chi_Q(\prod_{j=1}^k P_j^2)} = 1 -\frac 1{\#\Hgev} \#\{Q\in \Hgev: \exists P_j\mid Q\}$$
Replacing the set of square-free $Q$ by arbitrary monic $Q$ of
degree $2g+2$ gives
$$  \#\{Q\in \Hgev: \exists P_j\mid Q\} \leq
\#\{\deg Q=2g+2 :\exists P_j\mid Q\} \leq\sum_{j=1}^k\frac
{q^{2g+2}}{||P_j||}
$$
so that recalling $\Hgev = (1-\frac 1q) q^{2g+2}$, we have
$$ 1-   \frac 1{(1-\frac 1q) }\sum_{j=1}^k\frac{1}{||P_j||} \leq \ave{\chi_Q(\prod_{j=1}^k P_j^2)}  \leq 1
$$
Thus
$$
\ave{\chi_Q(\prod_{j=1}^k P_j^2)} = 1 +O\left( \sum_{j=1}^k
\frac1{||P_j||} \right)
$$
as claimed.
\end{proof}

For a polynomial $Q\in \FF_q[x]$ of positive degree, set
$$
\eta(Q) = \sum_{P\mid Q} \frac 1{||P||}
$$
the sum being over all monic irreducible (prime) polynomials dividing
$Q$.
\begin{lem}\label{lem:appendix}
The mean values of  $\eta$ and $\eta^2$
are uniformly  bounded as $g\to \infty$:
$$\ave{\eta}\leq 1, \qquad \ave{\eta^ 2} \leq \frac 1{(1-\frac 1q)^3 }
$$
\end{lem}
\begin{proof}
We consider the first moment:
We have
\begin{equation*}\begin{split}
\ave{\eta(Q)} &= \frac 1{\#\Hgev} \sum_{Q\in \Hgev} \sum_{P\mid Q} \frac
1{||P||} \\
&= \frac 1{\#\Hgev} \sum_{\deg P \leq 2g+2} \frac 1{||P||}
\#\{ Q\in \Hgev: P\mid Q \}
\end{split}\end{equation*}
We bound the number of
square-free $Q$ divisible by $P$ by the number of all $Q$ of degree
$2g+2$ divisible by $P$, which is $q^{2g+2}/||P||$, to find
\begin{equation*}\begin{split}
\ave{\eta(Q)}& \leq
\frac 1{(1-\frac 1q)q^{2g+2}} \sum_{\deg P \leq 2g+2} \frac 1{||P||}
\#\{ \deg Q=2g+2: P\mid Q \} \\
&\leq \frac 1{(1-\frac 1q) q^{2g+2}}
\sum_{\deg P \leq 2g+2}\frac{q^{2g+2}}{||P||^2}
\leq \frac 1{1-q^{-1}}\sum_{f} \frac 1{||f||^2} =  1
\end{split}\end{equation*}
(the last sum is over all monic polynomials) proving that
$\ave{\eta(Q)}$ is uniformly bounded.


For the second moment of $\eta$, we have
\begin{multline*}
\ave{\eta^2}=
\frac 1{\#\Hgev} \sum_{Q\in \Hgev} \left( \sum_{P\mid Q}
\frac{1}{||P||} \right)^2 \\
=\frac 1{\#\Hgev}  \sum_{\deg P_1,\deg P_2 \leq 2g+2}
    \frac 1{||P_1|| \cdot ||P_2||} \#\{Q\in \Hgev: P_1 \mid Q,
    P_2\mid Q\}
\end{multline*}
For squarefree $Q$, if two primes $P_1\mid Q$ and $P_2 \mid Q$ then necessarily
$P_1 \neq P_2$ and then $Q$ is divisible by both
iff it is divisible by their product, hence
\begin{equation*}
  \begin{split}
\#\{Q\in \Hgev: P_1\mid Q, P_2\mid Q \}
&=\#\{Q\in \Hgev: P_1P_2\mid Q\} \\
&\leq \#\{Q: \deg Q=2g+2, P_1P_2 \mid Q\}\\
& = \begin{cases}
\frac{q^{2g+2}}{||P_1 P_2||} ,& \deg(P_1P_2) \leq 2g+2\\0,&\mbox{otherwise}
\end{cases}
    \end{split}
\end{equation*}
and hence the contribution of such pairs is bounded by
$$
\frac 1{(1-\frac 1q)q^{2g+2}}
\sum_{P_1}\sum_{P_2}\frac{q^{2g+2}}{||P_1||^2||P_2||^2} \leq
\frac 1{(1-\frac 1q)} \left(\sum_{f} \frac 1{||f||^2} \right)^2 =
\frac 1{(1-\frac 1q)^3 }
$$
Thus we see $\ave{\eta^2} \leq (1-\frac 1q)^{-3}$
which is again uniformly bounded.
\end{proof}

\section{Beurling-Selberg functions}
Let $\mathcal I=[-\beta/2,\beta/2]$ be an interval, symmetric about
the origin, of length $0<\beta<1$, and $K\geq 1$ an integer.
Beurling-Selberg polynomials $I^{\pm}_K$ are trigonometric polynomials
approximating the indicator function $\mathbf 1_{\mathcal I}$
satisfying (see the beautiful exposition
in \cite[Chapter 1.2]{Montgomery ten}):
\begin{itemize}
\item $I^{\pm}_K$ are trigonometric polynomials of degree $\leq K$

\item Monotonicity:
\begin{equation}\label{beurling monotonicity}
I^-_K \leq \mathbf 1_{\mathcal I} \leq I^+_K
\end{equation}

\item
The integral of $I^{\pm}_K$ is close to the length of the
interval:
\begin{equation}\label{beurling2}
\int_0^1 I^{\pm}_K(x)dx =
\int_0^1 \mathbf 1_{\mathcal I}(x)dx \pm \frac 1{K+1}
\end{equation}

\item $I^{\pm}_K(x)$ are even\footnote{This is because we take the
interval ${\mathcal I} =[-\beta/2,\beta/2]$ which is symmetric about
the origin.}.
\end{itemize}

As a consequence of \eqref{beurling2},
the non-zero Fourier coefficients of $I^{\pm}_K$ satisfy
\begin{equation}\label{Fourier coeffs of I a}
\left| \^I^{\pm}_K(k) - \^{\mathbf 1}_{\mathcal I}(k) \right|
\leq \frac 1{K+1}
\end{equation}
and in particular
\begin{equation}\label{Fourier coeffs of I}
| \^I^{\pm}_K(k)| \leq \frac 1{K+1} + \min \left( \beta,\frac \pi
{|k|} \right) , \qquad  0<|k| \leq K
\end{equation}

\begin{prop}\label{prop sum of I(n)}
Let  $\mathcal I=[-\beta/2,\beta/2]$ be an interval and $K\geq 1$ an
integer so that $K\beta>1$.
Then
\begin{equation} \label{eq sum I(n)}
\sum_{n\geq 1} \^I_K^\pm(2n) = O(1)
\end{equation}
\begin{equation}\label{eq sum squares of I(n)}
\sum_{n \geq 1} n \^I_K^\pm(n)^2 =  \frac 1{2\pi^2} \log K\beta +O(1)
\end{equation}
where the implied constants are independent of $K$ and $\beta$.
\end{prop}
\begin{proof}
To bound the sum \eqref{eq sum I(n)},
we may use \eqref{Fourier coeffs of I a} to write
$$ \^I_K^\pm(2n) = \frac{\sin 2\pi n\beta}{2\pi n} + O(\frac 1K)
$$
and hence
$$
\sum_{n \geq 1} \^I_K^\pm(2n)= \sum_{1\leq n\leq K/2} \frac{\sin
  2\pi n\beta}{2\pi n}  +O(1)
$$

We treat separately the range $n<1/\beta$ and
$1/\beta<n<K$.
To bound the sum over $n<1/\beta$, use $\sin 2\pi n\beta \ll
n\beta$ and hence
$$ \sum_{1\leq n <1/\beta} \frac{\sin  2\pi n\beta}{2\pi n}
\ll \sum_{1\leq n <1/\beta} \frac{n\beta}{n} =O(1)
$$

For the sum on $n>1/\beta$, we apply summation by parts. The partial
sums of $\sin 2\pi n\beta$ are
\begin{equation}\label{sum by parts}
\sum_{n=1}^N \sin 2\pi n \beta =
\frac{ \cos \pi  \beta-\cos (2 N+1) \pi  \beta}{2 \sin \pi
  \beta}=O(\frac 1\beta)
\end{equation}
Therefore
$$
 \sum_{1/\beta<n<K/2} \frac{\sin  2\pi n\beta}{2\pi n} \ll
\frac 1{\beta K}+1+ \frac  1\beta \int_{1/\beta}^K \frac 1{t^2} dt =O(1)
$$
and hence $\sum_{n \geq 1} \^I_K^\pm(2n)=O(1)$.


To prove \eqref{eq sum squares of I(n)}, we use
\eqref{Fourier coeffs of I a} to write
$$
\sum_{n>0} n \^I_K^\pm(n)^2  =\frac 1{\pi^2} \sum_{n\leq K} \frac{(\sin \pi  n\beta )^2}{n}
+  O(1)
$$
We split the sum into two parts: The sum over $1\leq n \leq 1/\beta$,
where we use $|\sin \pi n \beta|\ll n\beta$ to see that it gives a
bounded contribution, and the sum over $1/\beta<n\leq K$, where we
use $\sin(y)^2 = \frac 12 (1-\cos(2y))$  to get
\begin{equation*}
\begin{split}
\sum_{n>0} n \^I_K^\pm(n)^2 &=
\frac 1{2\pi^2} \sum_{\frac 1\beta<n\leq K} \frac 1n
-\frac 1{2\pi^2} \sum_{\frac 1\beta< n\leq K} \frac{\cos 2\pi n\beta}{n}
+O(1) \\
&= \frac 1{2\pi^2} \log K\beta
-\frac 1{2\pi^2} \sum_{\frac 1\beta< n\leq K} \frac{\cos 2\pi n\beta}{n}
+O(1)
\end{split}
\end{equation*}
To bound $\sum_{\frac 1\beta< n\leq K} \frac{\cos 2\pi n\beta}{n}$,
apply summation by parts using
$$ \sum_{1\leq n\leq N} \cos 2\pi n\beta = \frac{\sin(2N+1)\pi \beta -
\sin \pi \beta}{2\sin \pi \beta } \ll \frac 1\beta, \quad 0<\beta <1
$$
to find that it gives a bounded contribution. Hence
$$
\sum_{n>0} n \^I_K^\pm(n)^2 =  \frac 1{2\pi^2}\log K\beta +O(1)
$$
as claimed.
\end{proof}

\section{Counting functions}
Let $\chi$ be a primitive Dirichlet character.
We denote by $N_{\mathcal I}(\chi)$ the number of angles
$\theta_{j,\chi}$ of the L-function $\mathcal L^*(u,\chi)$ 
(see \eqref{angles of L})  in the interval $\mathcal I=[-\beta/2,\beta/2]$.
Define $S_{\mathcal I}(\chi)$ by 
$$N_{\mathcal I}(\chi)=2g|\mathcal I| +\frac{2}{\pi} \arg
(1-\frac{e^{i\pi |\mathcal I|}}{\sqrt{q}})+ S_{\mathcal I}(\chi)$$

Set
$$N^\pm_{K}(\chi) = \sum_{j=1}^{D}
I^{\pm}_{K}(\theta_{j,\chi})
$$
Here $K$ will depend on $\deg Q$. This will be our approximation to
the counting function $N_{\mathcal I}(\chi)$. 
Then by virtue of \eqref{beurling monotonicity},
\begin{equation}\label{sandwich 1}
N^-_K(\chi) \leq N_{\mathcal I}(\chi) \leq N^+_K(\chi)
\end{equation}

Using the explicit formula \eqref{explicit formula}, we find
\begin{equation}\label{NSepsilon}
N^\pm_{K}(\chi) =: D(\beta\pm \frac 1{K+1}) +\lambda_\infty(\chi)
 \frac 1{\pi i} \int_0^1 I^{\pm}_K(\theta)
\frac{d}{d\theta} \log (1-\frac{e^{2\pi i\theta}}{\sqrt{q}}) d\theta
+ S^\pm_{K}(\chi)
\end{equation}
where $S^\pm_{K}(\chi)$ is
\begin{equation}\label{expr for Sepsilon}
S^\pm_{K}(\chi) :=
-\sum_{\deg f\leq K}  \^I^{\pm}_{K}(\deg f)
\frac{\Lambda(f)}{||f||^{1/2} }\left\{ \chi(f) +
\overline{\chi(f)} \right\}
\end{equation}
the sum taken over all prime powers  $f\in \FF_q[x]$ (of degree
$\leq K$).

Note that since $|| \mathbf 1_{\mathcal I} -I^{\pm}_K
||_{L^1} =\frac 1{K+1}$, we have
\begin{equation}\label{even contr}
\begin{split}
\frac 1{\pi i} \int_0^1 I^{\pm}_K(\theta)
\frac{d}{d\theta} \log (1-\frac{e^{2\pi i\theta}}{\sqrt{q}})
d\theta  &=
\frac 1{\pi i} \int_{-\beta/2}^{\beta/2}
\frac{d}{d\theta} \log (1-\frac{e^{2\pi i\theta}}{\sqrt{q}})
d\theta + O(\frac 1K)\\
& = \frac {2}{\pi} \arg (1-\frac{e^{i\pi \beta}}{\sqrt{q}})
 +O(\frac 1K)
\end{split}
\end{equation}

\subsection{Quadratic characters}
For the case at hand, of quadratic characters, we write 
$N_{\mathcal I}(Q)$ for $N_{\mathcal I}(\chi_Q)$, with similar meaning
for $S_{\mathcal I}(Q)$, $N_K^\pm(Q)$ and $S_K^\pm(Q)$. 
We have
\begin{equation}\label{expression for S(Q)}
 S_K^\pm(Q):= S^\pm_{K}(\chi_Q) =
-2\sum_{\deg f\leq K}  \^I^{\pm}_{K}(\deg f)
\frac{\Lambda(f)}{||f||^{1/2} }\chi_Q(f)
\end{equation}

We may now deduce that the zeros are uniformly distributed:
\begin{prop}\label{prop:unif dist}
Every fixed (symmetric) interval $\mathcal I=[-\beta/2,\beta/2]$ contains
asymptotically $2g|\mathcal I|$ angles $\theta_{j,Q}$, in fact
$$ N_{\mathcal I}(Q) = 2g|\mathcal I| + O(\frac g{\log g})$$
\end{prop}
\begin{proof}
Indeed from \eqref{sandwich 1}  it suffices to show that for the
smooth counting functions $N^\pm_K(\chi_Q)$ we have
$$ N^\pm_K(\chi_Q) = 2g|\mathcal I| + O(\frac g{\log g})$$
Now from \eqref{NSepsilon}, \eqref{even contr} it follows that
$$
N^\pm_K(\chi_Q) = 2g|\mathcal I| +O(\frac gK) +O(1) + |S_K^\pm(Q)|
$$
To bound $S_K^\pm(Q)$, use \eqref{expression for S(Q)}
and \eqref{Fourier coeffs of I} in the form
$\^I_K^\pm(\deg f)\Lambda(f) = O(1)$ to deduce that
$$
S_K^\pm(Q) \ll \sum_{\deg f\leq K} \frac 1{\sqrt{||f||}} \ll q^{K/2}
$$
and hence
$$
\left| N^\pm_K(\chi_Q) - 2g|\mathcal I| \right| \ll \frac gK +q^{K/2}
$$
Taking $K\approx \log_q g-\log_q\log g$ gives the result.
\end{proof}

\section{Expected value}
We first bound the expected value of $S_{\mathcal I}$:
\begin{prop}\label{prop: expectation}
Assume that either the interval $\mathcal I=[-\beta/2,\beta/2]$ is
fixed or that it shrinks to zero with $g\to \infty$ in such a way that
$g\beta\to \infty$. Then
$$\ave{S_{\mathcal I}} =O(1)$$
\end{prop}
\begin{proof}
Using \eqref{sandwich 1}, \eqref{NSepsilon} and \eqref{even contr},
we find that for any $K$,
$$\ave{S_K^-}\leq \ave{S}+O( \frac {g}{K})\leq \ave{S_K^+}$$
Taking $K\approx g/100$ gives the remainder term above is bounded. So
it remains to bound the expected value of $S_K^{\pm}$ for such $K$.

Recall that $S_K^{\pm}$ is a sum over prime powers. We separate out
the contribution of even powers, which is not oscillatory, from that
of the odd powers:
$$
S_K^{\pm} = \mbox{ even} + \mbox{ odd}
$$
We claim that the even powers give
\begin{equation}\label{eq: expression for even}
\mbox{ even} = -2 \sum_{n \geq 1} \^I_K^\pm(2n)  +O\left(\eta(Q) \right)
\end{equation}
where
$$ \eta(Q) = \sum_{P\mid Q} \frac 1{||P||}$$
the sum over prime divisors of $Q$.

To see \eqref{eq: expression for even}, note that  for an even power
of a prime, say $f=g^2$, we have
$\chi_Q(f)=1$ if $\gcd(g,Q)=1$ and $0$ otherwise.
Writing the even powers of a prime as $f=g^2$, and noting that
$\Lambda(f)=\Lambda(g)$, we have
\begin{equation*}
  \begin{split}
    \mbox{even } &= -2\sum_{\gcd(g,Q)=1} \frac{\^I_K^\pm(2\deg g)
      \Lambda(g)}{||g||} \\
 &= -2\sum_{n\geq 1} \frac{\^I_K^\pm(2n) }{q^n} \sum_{\deg g=n}
    \Lambda(g)  + O\left( \sum_{P\mid Q} \frac{1}{||P||}
    \right)
  \end{split}
\end{equation*}
where the remainder term is a sum over all prime divisors of $Q$.
By the prime number theorem, $\sum_{\deg g=n} \Lambda(g)=q^n$ and
hence
\begin{equation}
-2\sum_{\gcd(g,Q)=1} \frac{\^I_K^\pm(2\deg g) \Lambda(g)}{||g||}=
 -2\sum_{n \geq 1} \^I_K^\pm(2n)
\end{equation}
proving \eqref{eq: expression for even}.

It now follows that expected value of the even powers is bounded:
Indeed, the sum $\sum_{n \geq 1} \^I_K^\pm(2n)$
is  bounded by Proposition~\ref{prop sum of I(n)} (note that our
choice $K\approx g/100$ and the condition $g\beta\to\infty$ guarantees
$K\beta\to \infty$, hence Proposition~\ref{prop sum of I(n)} is
applicable). As for the term $\eta(Q)=\sum_{P\mid Q} \frac{1}{||P||}$,
it is not bounded individually, but its mean is bounded by
Lemma~\ref{lem:appendix}.

The expected value of the odd powers is
$$
\ave{\mbox{ odd}}  = -2\sum_{\deg f \mbox{ odd}} \frac{\^I_K^\pm(\deg f)
\Lambda(f)}{\sqrt{||f||}} \ave{\chi_Q(f)}
$$
To estimate the expected value of the odd powers,
we use Lemma~\ref{lem:average chi} and  \eqref{Fourier coeffs of I}
in the form  $\^I_K^\pm(\deg f)\Lambda(f) = O(1)$  to find
$$
\ave{\mbox{ odd}} \ll \sum_{\deg f\leq K}\frac 1{\sqrt{||f||}}
\frac{2^{\deg f}}{q^{g+1}}  \ll \frac{(2\sqrt{q})^K}{q^{g+1}}
$$
which for $K\approx g/100$ is bounded.
\end{proof}
Hence we see that
\begin{equation}
\ave{ \frac{S_{\mathcal I}}{\sqrt{\frac 2{\pi^2}\log(g\beta)}}} \to 0, \qquad
g\to \infty
\end{equation}


\section{A sum over primes}

Consider the sum over primes
$$
T_K^\pm(Q):=- 2\sum_P \frac{\^I_K^\pm(\deg P) \deg P} {\sqrt{||P||}}
\chi_Q(P)
$$
This will be our approximation to $S_{\mathcal I}$. 
From now on assume that
$$
K \approx \frac{g}{\log\log(g\beta)}
$$
which will guarantee $\log K\beta \sim \log g\beta$ and $K=o(g)$.


\begin{thm}\label{prop approximation by prime sum}
Assume that $g\to \infty$ and either $0<\beta<1$ is fixed or $\beta\to 0$
while $\beta g\to \infty$. Take $K\approx g/\log\log(g\beta)$. Then

i) $$
\ave{|T_K^\pm|^2} \sim  \frac 2{\pi^2}  \log \beta g
$$

  ii)
  \begin{equation}\label{difference of T's}
   \ave{|T_K^+-T_K^-|^2}=O(1)
  \end{equation}

  iii)
  \begin{equation}\label{difference of S and T}
  \ave{|S_K^\pm - T_K^\pm|^2} =O(1)
  \end{equation}
  \end{thm}
  The rest of this section is devoted to the proof of
  Theorem~\ref{prop approximation by prime sum}.




  \subsection{Computing $\ave{\left(T^\pm_{K}\right)^2}$}\label{subsec:7.1}
  We have
  $$
  \ave{\left(T^\pm_{K}\right)^2}
  =4\sum_{P_1,P_2} \^I^{\pm}_{K}(\deg P_1)\^I^{\pm}_{K}(\deg P_2)
  \frac{\deg P_1 \deg P_2}{\sqrt{||P_1|| ||P_2||}}
  \ave{ \chi_Q(P_1 P_2)}
  $$
  The sum is  over $\deg P_1,\deg P_2\leq K <g$.
  Consider the contribution of pairs such that $P_1P_2$ is not a perfect
  square (the ``off-diagonal pairs'').
  We may use Lemma~\ref{lem:average chi} to bound their contribution by
  $$\ll \frac 1{q^{g+1} }  \left( \sum_{\deg P \leq K}
   \frac{ |\^I^{\pm}_{K}(\deg P)|\deg P 2^{\deg P}} {\sqrt{|P|}}  \right)^2
  $$
  Using \eqref{Fourier coeffs of I} in the form $|\^I^{\pm}_{K}(k)|\ll
  1/|k|$ gives that  the inner sum is bounded by
  $$ \ll  \sum_{\deg P \leq K} \frac {2^{\deg P}}{\sqrt{|P|}} \frac
  {\deg P}{\deg P} \ll (2\sqrt{q})^{K}
  $$
  Hence the off-diagonal contribution is  bounded by
  $$ \ll \frac {(4q)^K}{q^{g+1} }  $$
  which is negligible since we take $K=o(g)$.

  Consider the contribution of pairs such that $P_1 \cdot P_2$ is a
  square. Since $P_1$ and $P_2$ are primes, this forces $P_1=P_2$.
These contribute
\begin{multline}\label{diagonal}
4\sum_P \frac{(\deg P)^2}{||P||} \^I^{\pm}_{K}(\deg P)^2\ave{\chi_Q(P)^2} \\
=4\sum_P \frac{(\deg P)^2}{||P||} \^I^{\pm}_{K}(\deg P)^2
+ O\left( \sum_P \frac{(\deg P)^2}{||P||^2} \^I^{\pm}_{K}(\deg P)^2  \right)
\end{multline}
by Lemma~\ref{lem:ave chi^2}.

Using the prime number theorem $\#\{P:\deg P=n\}= q^n/n+O(q^{n/2})$ gives
\begin{equation*}
\begin{split}
  4\sum_P \frac{(\deg P)^2}{||P||} \^I^{\pm}_{K}(\deg P)^2 &=
  4\sum_{1\leq n\leq K}  \left(n+O(\frac{n^2}{q^{n/2}})\right)
  \^I^{\pm}_K(n)^2 +O(1)\\
&= 4\sum_{1\leq n\leq K} n \^I^{\pm}_K(n)^2 +O(1)
\end{split}
\end{equation*}
By Proposition~\ref{prop sum of I(n)}  we find 
\begin{equation}\label{variance sum} 
4\sum_P \frac{(\deg P)^2}{||P||} \^I^{\pm}_{K}(\deg P)^2= 
\frac 2{\pi^2}\log K\beta +O(1)  
\end{equation}
(note that if $g\beta\to \infty$
then  $K\beta\approx g\beta/\log\log(g\beta)\to\infty$).
To bound the remainder term in \eqref{diagonal} use
\eqref{Fourier coeffs of I} in the form $\^I^{\pm}_{K}(\deg P) \deg P = O(1)$
to find that the sum is at most $\sum_P 1/||P||^2= O(1)$.
Therefore we find
$$
\ave{\left(T^\pm_{K}\right)^2} = \frac 2{\pi^2}\log (K\beta) +O(1)\;.
$$

\subsection{Bounding $\ave{ |T^+_{K} - T^-_{K}|^2 }$}
  Next we compute the variance of  the difference
  $ \ave{ \left| T^+_{K} -T^-_{K} \right|^2 }$.
  Arguing as above, one sees that the only terms which may significantly
  contribute to the average are again the diagonal terms
  $$
  \ave{\left| T^+_{K} - T^-_{K} \right|^2 } = 4\sum_{\deg P\leq K}
  \frac{(\deg P)^2}{||P||} \left( \^I^+_K(\deg P)-\^I^-_{K}(\deg P)
  \right)^2 \ave{\chi_Q(P)^2} +o(1)
  $$
  Since by \eqref{Fourier coeffs of I a}
  $$
  \left|  \^I^+_K(n) - \^I^-_K(n) \right| \leq \frac 2{K+1}
  $$
  we get
  $$
   \ave{\left| T^+_{K} - T^-_{K} \right|^2 }
   \ll \frac 1{K^2} \sum_{\deg P \leq K} \frac{(\deg P)^2}{||P||}
  $$
  Using the Prime Number Theorem, this is easily seen to be
  $O(1)$. Hence we find
  $$
   \ave{\left| T^+_{K} - T^-_{K} \right|^2 }
   =O(1)
  $$

  \subsection{Bounding $\ave{|S_K^\pm - T_K^\pm|^2}$}
  Next we  show that $\ave{|S_K^\pm - T_K^\pm|^2} =O(1)$.
  We have
  \begin{equation}
    \begin{split}
  S_K^\pm - T_K^\pm &= -2\sum_{f=P^j, j\geq 2}
  \^I^{\pm}_{K}(\deg f) \frac{\Lambda(f)}{||f||^{1/2} }\chi_Q(f)   \\
  & = \mbox{even} + \mbox{odd}
    \end{split}
  \end{equation}
  where the term ``even'' is a sum over the even powers of primes, and
  ``odd'' is the sum over odd powers of primes where the exponent is at
  least $3$.  We will show that the second moments of both the odd and
  even terms are bounded.

  We first argue that the second moment of the even powers contribute a
  bounded amount.
  As we saw in the proof of Proposition~\ref{prop: expectation}, see
  \eqref{eq: expression for even}, we have
  $$\mbox{even} \ll 1+ \sum_{P\mid Q} \frac 1{||P||}$$
  the sum being over all prime divisors of $Q$.
  This is not bounded individually, but its second moment is
  bounded by Lemma~\ref{lem:appendix}.

  It remains to bound the contribution of the odd powers.
  We have
  $$
  \ave{|\mbox{odd}|^2} = 4 \sum_{ f_1, f_2}
  \^I^{\pm}_{K}(\deg f_1) \^I^{\pm}_{K}(\deg
  f_2)\frac{\Lambda(f_1)\Lambda(f_2)}{||f_1 f_2||^{1/2} }
  \ave{\chi_Q(f_1 f_2)}
  $$
  where the sum is over odd higher prime powers, that is over $f=P^j$
  with $j\geq 3$ and odd.

  The pairs where $f_1\cdot f_2$ is not a square contribute $o(1)$ by the
  same argument as above.
  Consider the contribution of pairs such that $f_1 \cdot f_2$ is a
  square.
  If $f_1$ and $f_2$ are {\em odd} higher prime powers but $f_1 \cdot f_2$ is a
  square, then necessarily  $f_1=P^r$, $f_2=P^s$ with $P$ prime, $r,s \geq 2$,
  (and $r=s \mod 2$).
  Necessarily then $r+s\geq 4$.
  The contribution of such pairs can be bounded, using
  \eqref{Fourier coeffs of I} in the form
  $\^I^\pm_K(\deg f) \Lambda(f)=O(1)$, by
  $$
  \sum_P \sum_{ r+s\geq 4} \frac 1{||P||^{(r+s)/2}}\ll
  \sum_P \sum_{j\geq 4} \frac{j}{||P||^{j/2}}
  \ll \sum_P \frac 1{||P||^2} = O(1)
  $$
  Hence  $\ave{|\mbox{odd}|^2}=O(1)$ and therefore
  $$\ave{|S_K^\pm - T_K^\pm|^2} =O(1)$$

\section{Higher moments of  $T^\pm_K$}
In this section we show that all moments of $T^\pm_K$ are Gaussian.
\begin{thm}\label{thm: higher moments of T}
Assume the setting of Theorem~\ref{prop approximation by prime sum}
and let $r\geq2$. Then
$$
\left|\ave{(T^\pm_K)^{2r-1}}\right|
=o(1) $$
and
$$
\ave{(T^\pm_K)^{2r}}=\frac{(2r)!}{r!\pi^{2r}}\log^r (\beta
K)+O\left(\log^{r-1} (\beta  K)\right)
$$
\end{thm}
\begin{proof}
  For the odd moments, we have
$$
\ave{\left(T^\pm_{K}\right)^{2r-1}}
  =-2^{2r-1}\sum_{P_1,...,P_{2r-1}}\frac{\prod \^I^{\pm}_{K}(\deg P_j)\deg P_j}{\sqrt{||\prod P_j|| }}
  \ave{ \chi_Q(\prod P_j)}
    $$
    Since $\prod_j P_j$ cannot be a perfect square, we may apply
    lemma~\ref{lem:average chi} and obtain the bound
  $$\left|\ave{\left(T^\pm_{K}\right)^{2r-1}}\right|\ll \frac 1{q^{g+1} }  \left( \sum_{\deg P \leq K}
   \frac{ |\^I^{\pm}_{K}(\deg P)|\deg P 2^{\deg P}} {\sqrt{|P|}}
   \right)^{2r-1}
  $$
  As was already calculated in \S~\ref{subsec:7.1}, the inner sum is
  bounded by
  $$ \ll  \sum_{\deg P \leq K} \frac {2^{\deg P}}{\sqrt{|P|}} \frac
  {\deg P}{\deg P} \ll (2\sqrt{q})^{K}
  $$
  Hence
  $$ \left|\ave{\left(T^\pm_{K}\right)^{2r-1}}\right|\ll \frac {(2\sqrt q)^{(2r-1)K}}{q^{g+1} }$$
which vanishes assuming $K\approx g/\log\log(g\beta)$.

To compute the even moments, write
$$
\ave{\left(T^\pm_{K}\right)^{2r}}=2^{2r}(T_{sq}^{2r}+T_{nsq}^{2r})
$$
where both $T_{sq}^{2r}$ and $T_{nsq}^{2r}$ have the form
$$
\sum_{P_1,...,P_{2r}}\frac{\prod \^I^{\pm}_{K}(\deg P_j)\deg
  P_j}{\sqrt{||\prod P_j|| }}   \ave{ \prod\chi_Q(P_j)}
$$
where $T_{sq}^{2r}$ is the sum over prime $2r$-tuples $\{P_j\}$ for
which $\prod_{j=1}^{2r}P_j$ is a perfect square, and $T_{nsq}^{2r}$
contains the  remaining (off-diagonal) terms.

The term  $T_{nsq}^{2r}$ can be bounded as was done for the odd moments:
$$
T_{nsq}^{2r}\ll  \frac 1{q^{g+1} }  \left( \sum_{\deg P \leq K}
   \frac{ |\^I^{\pm}_{K}(\deg P)|\deg P 2^{\deg P}} {\sqrt{|P|}}
   \right)^{2r}\ll \frac {(2\sqrt q)^{2rK}}{q^{g+1} }
$$

Now
$$
T_{sq}^{2r}=\sum_{P_1\cdot \dots \cdot P_{2r}=\Box}\frac{\prod \^I^{\pm}_{K}(\deg
  P_j)\deg P_j}{\sqrt{||\prod P_j||   }}\ave{ \prod\chi_Q(P_j)}
$$
the sum taken over only those primes for which $\prod P_j$ is a
square, which implies all $P_j$ appear in equal pairs in each
summand. Note that in particular all summands are positive. By lemma
~\ref{lem:ave chi^2} we may replace $\ave{ \prod\chi_Q(P_j)}$ with
$1$ by introducing an error of $O\left(\sum_j1/||P_j||\right)$.

The total error produced by this substitution is, keeping in mind that
the primes $P_1,\dots P_{2r}$ must come in identical pairs, bounded by
$$ \sum_{j=1}^r \sum_{P_1,...,P_r} \frac{\prod_{k=1}^r\^I^{\pm}_{K}(\deg P_k)^2 (\deg
P_k)^2}{||P_j||^2\prod_{k\neq j}||P_k||}\ll$$
$$\ll
\sum_{P_2,...,P_r} \frac{\prod_{k=2}^r\^I^{\pm}_{K}(\deg P_k)^2
(\deg
P_k)^2}{\prod_{k=2}^r||P_k||}\sum_{P_1}\frac{\^I^{\pm}_{K}(\deg
P_1)^2 (\deg P_1)^2}{||P_1||^2}
$$ 
The inner sum is bounded, and hence the total error introduced is
$$
\ll\sum_{P_2,...,P_r} \frac{\prod_{k=2}^r\^I^{\pm}_{K}(\deg P_k)^2
(\deg P_k)^2}{\prod_{k=2}^r||P_k||}\ll \left( \log(\beta K)\right)^{r-1}
$$ 
by \eqref{variance sum}.

So far we showed that
$$
T_{sq}^{2r}=\sum_{P_1\cdot \dots \cdot P_{2r}=\Box } \frac{\prod
\^I^{\pm}_{K}(\deg P_j)\deg P_j}{\sqrt{||\prod
P_j||}}+O(\log^{r-1}(\beta K))
$$
Now we show that pairs of equal $P_j$ in
$$\sum_{P_1\cdot \dots \cdot P_{2r}=\Box}
\frac{\prod \^I^{\pm}_{K}(\deg P_j)\deg P_j}{\sqrt{||\prod P_j||}}
$$
can be taken all distinct, for the remaining terms are bounded by
\begin{multline*}
\ll\sum_{P_1=P_2=P_3=P_4}\frac{\^I^{\pm}_{K}(\deg P_1)^4\deg^4
P_1}{||P_1||^2}\sum_{\prod_{j=5}^{2r}P_j=\Box}\frac{\prod
\^I^{\pm}_{K}(\deg P_j)\deg P_j}{\sqrt{||\prod
P_j||}}\\\ll\sum_{j=0}^\infty\frac{q^j}{j}\frac{j^4}{q^{2j}}
\log^{r-2}(\beta K)\ll\log^{r-2}(\beta K)
\end{multline*}
Finally, the sum over distinct pairs is
$$
\frac{(2r)!}{r!2^r}\sum_{P_1,...,P_r\mbox{ distinct}}
\frac{\prod \^I^{\pm}_{K}(\deg P_j)^2\deg^2 P_j}{||\prod P_j||}
$$
Now we remove the restriction that $P_1,\dots P_r$ are
distinct, introducing (again) an error of $O(\log^{r-2}(\beta K))$,
and obtain
$$ T_{sq}^{2r}=\frac{(2r)!}{r!2^r}\left(\sum_P\frac{\^I^{\pm}_{K}(\deg
P)^2\deg^2 P}{||P||}\right)^r + O(\log^{r-1}(\beta K))
$$

Summarizing all said above, and using \eqref{variance sum} 
yields
$$
T_{sq}^{2r}
=\frac{(2r)!}{r!\pi^{2r}2^{2r}}\log^r(\beta K) +O(\log^{r-1}(\beta K))
$$
and
$$
\ave{(T^\pm_K)^{2r}}=\frac{(2r)!}{r!\pi^{2r}}\log^r (\beta K)
+O\left(\log^{r-1} (\beta  K)\right)
$$
as claimed.
\end{proof}


\begin{cor}\label{cor: moments of T}
Under the assumption of Theorem~\ref{prop approximation by prime sum},
$T^\pm_K/\sqrt{\frac 2{\pi^2}  \log  g\beta}$
has a standard Gaussian limiting distribution.
\end{cor}
Indeed, the main-term expressions for the moments of $T^\pm_K$ imply
all moments of $T^\pm_K/\sqrt{\frac 2{\pi^2}  \log g\beta}$ are
  asymptotic to standard Gaussian moments, where the odd moments
vanish and the even moments are
$$ \frac{1}{\sqrt{2\pi}} \intinf x^{2r} e^{-x^2/2}dx =
1\cdot 3\cdot \dots \cdot (2r-1)=\frac{(2r)!}{2^r r!}
$$


\section{Conclusion}

In this section we prove the claim \eqref{CLT for S} in our
introduction.
Recall that we wrote
$$N_{\mathcal I}(Q)=2g|\mathcal I| +\frac{2}{\pi} \arg
(1-\frac{e^{i\pi |\mathcal I|}}{\sqrt{q}})+ S_{\mathcal I}(Q)$$
and thus \eqref{CLT for S}  is equivalent to:
\begin{thm}\label{main thm}
Assume either that the interval $\mathcal I=[-\beta/2,\beta/2]$ is
fixed, or that its length $\beta$ shrinks to zero while
$g\beta\to\infty$. Then
\begin{equation*}
\ave{|S_{\mathcal I}|^2} \sim \frac 2{\pi^2}  \log g\beta
\end{equation*}
and $S_{\mathcal I}/\sqrt{\frac 2{\pi^2}\log \beta g}$
has a standard Gaussian  distribution.
\end{thm}
To prove this, it suffices to show that the second moment of the difference
  $S_{\mathcal I}-T^\pm_K$ is negligible relative to $ \log(g\beta)$:
  \begin{prop}\label{thm second moment}
  Assume that $K\approx g/\log\log g\beta$, and that either $\beta$ is
  fixed or $\beta\to 0$ while $g\beta\to\infty$. Then
  \begin{equation}
   \ave{ \left|\frac{S_{\mathcal I}-T^\pm_K}{\sqrt{ \frac 2{\pi^2}  \log g\beta }}
  \right|^2} \to 0
  \end{equation}
  \end{prop}
  Indeed, due to Proposition~\ref{thm second moment},
  the second moment of $S_{\mathcal I}$ is close to that of $T_K^\pm$ and
  and  the distribution of
  $S_{\mathcal I}/\sqrt{\frac 2{\pi^2}\log \beta g}$
  coincides with that of $T_K^\pm/\sqrt{\frac 2{\pi^2} \log (g\beta)}$,
  that is by Corollary~\ref{cor: moments of T}
we find that $S_{\mathcal I}/\sqrt{\frac 2{\pi^2}\log \beta g}$
  has a standard Gaussian distribution. Thus we will have proved
  Theorem~\ref{main thm} once we establish Proposition~\ref{thm second moment}.

  \subsection{Proof of Proposition~\ref{thm second moment}}
  Assume that $K\approx g/\log\log(g\beta)$. Then it suffices to show
  \begin{equation} \label{eq: S - T_K}
      \ave{ |S_{\mathcal I}-T^\pm_K|^2} \ll (\frac gK)^2 \;.
  \end{equation}

  We first show
    \begin{equation}\label{var of S-S+}
  \ave{|S_{\mathcal I} -S^\pm_K|^2} \ll (\frac{g}{K})^2     \;.
    \end{equation}
  By \eqref{sandwich 1}, we have
  $$
  S^-_K \leq S_{\mathcal I}+O(\frac gK) \leq S^+_K \;.
  $$
  Hence
  $$ 0\leq S_{\mathcal I} -S^-_K +O(\frac gK) \leq S^+_K - S^-_K \;.
  $$
  Since we are dealing now with positive quantities, we may take
  absolute values and get
  $$ |  S_{\mathcal I}-S^-_K +O(\frac gK)| \leq |S^+_K - S^-_K|
  $$
  and applying the triangle inequality gives
  $$ |S_{\mathcal I}-S^-_K| \leq |S^+_K-S^-_K| +O( \frac gK) \;,
  $$
  hence
  $$  |S_{\mathcal I}-S^-_K|^2 \leq 2|S^+_K-S^-_K|^2 +O( (\frac gK)^2) \;.
  $$
  Taking expected values we get
  \begin{equation}\label{eq:temp bound}
   \ave{|S_{\mathcal I}-S^-_K|^2} \leq 2\ave{|S^+_K- S^-_K|^2} +
  O\left( (\frac gK)^2 \right) \;.
  \end{equation}

  To bound $\ave{|S_K^+-S_K^-|^2}$, use the triangle inequality to get
  $$
  |S_K^+-S_K^-| \leq  |S_K^+ - T_K^+| + |T_K^+-T_K^-| + |T_K^- -S_K^-|
  $$
  and hence
  $$
  |S_K^+-S_K^-|^2 \leq  3 \left( |S_K^+ - T_K^+|^2 + |T_K^+-T_K^-|^2 +
   |T_K^- -S_K^-|^2 \right) \;.
  $$
  Applying \eqref{difference of T's} and \eqref{difference of S and T} we find
  \begin{equation}\label{dif of S+ - S-}
   \ave{|S_K^+-S_K^-|^2} = O(1) \;.
  \end{equation}
  Inserting  \eqref{dif of S+ - S-} into \eqref{eq:temp bound} gives
  $$  \ave{|S_{\mathcal I}-S^-_K|^2} \ll ( \frac gK)^2$$
  and together with \eqref{dif of S+ - S-} we get
  $$  \ave{|S_{\mathcal I}-S^+_K|^2} \ll ( \frac gK)^2$$
  proving \eqref{var of S-S+}.

  To show \eqref{eq: S - T_K}, we use the triangle inequality to get
  $$ |S_{\mathcal I}-T^\pm_K| \leq |S_{\mathcal I}-S^\pm_K| + |S^\pm_K-T^\pm_K|$$
  hence
  $$
  \ave{|S_{\mathcal I}-T^\pm_K|^2} \leq 2\ave{|S_{\mathcal I}-S^\pm_K|^2} +
  2\ave{|S^\pm_K-T^\pm_K|^2}
  $$
  which is $O((\frac gK)^2 )$ by \eqref{var of S-S+} and
  \eqref{difference of S and T}.
  \qed

\end{document}